\documentclass[a4paper,10pt]{article}
\usepackage[dvips]{graphicx}
\usepackage{amsmath}
\usepackage{a4wide}
\usepackage{amssymb}
\usepackage{amsfonts}
\usepackage{epsfig}
\usepackage{color}

\newcommand{\al}{\mathfrak{g}}

\newcommand{\G}{\mathbb{G}}
\newcommand{\C}{\mathcal{C}}
\newcommand{\RIC}{RIC }
\newcommand{\LIC}{LIC }
\newcommand{\TC}{TC }

\newcommand{\Ad}{\mathrm{Ad}}

\title{Coordinated motion design on Lie groups}
\author{A. Sarlette, S. Bonnabel and R. Sepulchre\thanks{The authors are with the Department of Electrical Engineering and Computer Science, University of Li\`ege, Belgium. This paper presents research results of the Belgian Network DYSCO (Dynamical Systems, Control, and Optimization), funded by the Interuniversity Attraction Poles Programme, initiated by the Belgian State, Science Policy Office. The scientific responsibility rests with its authors. A. Sarlette is supported as an FNRS fellow (Belgian Fund for Scientific Research).}}

\begin{document}
\definecolor{mygray}{rgb}{0.6,0.6,0.6}

\maketitle

\begin{abstract}
The present paper proposes a unified geometric framework for coordinated motion on Lie groups. It first gives a general problem formulation and analyzes ensuing conditions for coordinated motion. Then, it introduces a precise method to design control laws in fully actuated and underactuated settings with simple integrator dynamics. It thereby shows that coordination can be studied in a systematic way once the Lie group geometry of the configuration space is well characterized. This allows among others to retrieve control laws in the literature for particular examples. A link with Brockett's double bracket flows is also made. The concepts are illustrated on $SO(3)$, $SE(2)$ and $SE(3)$.
\end{abstract}



\section{Introduction}

Recently, many efforts have been devoted to the design and analysis of control laws that coordinate swarms of identical autonomous agents --- see e.g. oscillator synchronization \cite{Sync,StrKur}, flocking mechanisms \cite{JADBABAIE,hendrickx1}, vehicle formations \cite{FaxMurray,DesaiKumar,OlfatiMurray,JandK,JandK2D,JandK3D}, spacecraft formations \cite{HadaeghSatSynch,VanDyke1,BeardOnSats,SO3book,MCINNES,Dimarogonas1,WRen1}, mechanical system networks  \cite{SHNOrControl,HNSyms,SujitThesis} and mobile sensor networks \cite{LUCA2,SPL2005,TAC2,LEONARDgen,DanOscillatory}.
For systems on vector spaces, so-called \emph{consensus algorithms} have been shown to be efficient and robust \cite{MOREAU,MOREAU2,ConsensusReview,olfati,JADBABAIE,TsitsiklisThesis}, and allow to address many relevant engineering issues and tasks \cite{MOREAU,FaxMurray,Tsitsiklis2}. However, in many of the above applications, the agents to coordinate evolve on nonlinear manifolds: oscillators evolve on the circle $S^1 \cong SO(2)$, satellite attitudes on $SO(3)$ and vehicles move in $SE(2)$ or $SE(3)$; these particular manifolds actually share the geometric structure of a \emph{Lie groups}. Coordination on nonlinear manifolds is inherently more difficult than on vector spaces. The goal of the present paper is to propose a unified geometric framework for coordinated motion on Lie groups, from a geometric definition of ``coordination'' to a purely geometric derivation of control laws for coordination like those proposed in \cite{SPL2005,TAC2,LUCA2,BaiWen,MY3,MY5}, in fully actuated and underactuated settings with simple integrator dynamics.\\

\paragraph*{Symmetries} The starting point for the developments in this paper is to assume \emph{invariance} or \emph{symmetry} in the behavior of the swarm of agents with respect to their absolute position on the Lie group: only \emph{relative} positions of the agents matter. In the 3-dimensional physical world, the laws governing interactions in a set of particles are invariant with respect to translations and rotations of the whole set as a rigid body. From this viewpoint, the invariance assumption comes down to assuming that there is no external influence acting on the agents. The symmetries of the system determine how to define meaningful quantities for the swarm, like ``relative positions'', and what the dynamics of the coupled agents can be. \emph{Coordinated motion} --- in short \emph{coordination}--- is defined as all situations where the relative positions of the agents are fixed. Feedback control laws that asymptotically enforce coordination must be designed on the basis of error measurements involving appropriately invariant quantities (e.g. \emph{relative} agent positions) only.\\

\paragraph*{Previous work} Results about synchronization (``reaching a common point'') and coordinated motion (``moving in an organized way'') on vector spaces are becoming well established \cite{TsitsiklisThesis,MOREAU,olfati,ConsensusReview}. Because a vector space can be identified with its tangent plane, both synchronization and coordinated motion can be seen as consensus problems on the same vector space: the former is a \emph{position} consensus while the latter is a \emph{velocity} consensus. In contrast when the configuration space is a Lie group, synchronization and coordinated motion are fundamentally different things. The geometric viewpoint for dynamical systems on Lie groups is a very well studied subject; see basic results in \cite{Jurdjevic,ArnoldMech} for simplified dynamics like those considered in the present paper, and \cite{MarsdenBook2,LeonardThesis,BulloThesis,ArnoldMech} for a geometric theory of \emph{mechanical} systems on Lie groups. General results for \emph{synchronization} on compact Lie groups are proposed in \cite{MY4}; see also \cite{MY4} for links to related examples in the literature. But to the best of the authors' knowledge, a unified geometric viewpoint for \emph{coordinated motion} --- in short \emph{coordination} --- on Lie groups is still lacking. Close to the present paper in its geometric flavor, \cite{Silvere1} builds invariant \emph{observers} for systems with Lie group symmetries; observer design can be seen as two-agent leader-follower synchronization on Lie groups.

In applications, the ubiquitous example of motion on Lie groups is a rigid body in $\mathbb{R}^n$. When translational motion is discarded, the configuration space reduces to the compact Lie group $SO(n)$; an element of $SO(n)$ can be represented by the $n \times n$ rotation matrix between a frame attached to the rigid body and a hypothetical fixed reference frame. The standard example of this type is satellite attitude control, where \emph{synchronization}, i.e. obtaining equal orientations, has recently attracted much attention \cite{HadaeghSatSynch,VanDyke1,SO3book,NorwayLeaderFollower,BeardOnSats,EADSDarwinControl,Sujit1,SujitThesis,Sujit3,MY3,Dimarogonas1,WRen1}, with and without external reference tracking; note that synchronization is a very special case of coordination. Considering rotations \emph{and translations}, the configuration space of an $n$-dimensional rigid body becomes the non-compact Lie group $SE(n) = \mathbb{R}^n \ltimes SO(n)$. Recently, coordination has been investigated on $SE(2)$ \cite{JandK2D,SPL2005,TAC2} and $SE(3)$ \cite{JandK3D,LUCA2,SHNOrControl,HNSyms} in the underactuated setting of \emph{steering control} where the linear velocity is fixed in the body's frame. Motion on $SE(n)$ with steering control is also directly linked to the evolution of a Serret-Frenet frame with curvature control, as explained in \cite{Jurdjevic}. Results taking into account the full mechanical dynamics for rigid body motion are more difficult to obtain --- see for instance applications of the framework of \cite{MarsdenBook2} for coordination on $SO(3)$ and $SE(3)$ in \cite{SujitThesis,Sujit3} and \cite{SHNOrControl,HNSyms} respectively. Considering simplified dynamics, as in the present paper, can be useful either to build a high-level planning controller or as a preliminary step towards an integrated mechanical controller, as illustrated for synchronization on $SO(3)$ in \cite{MY3} and \cite{MY5,MY6} respectively.\\

\paragraph*{Contributions} The main goal of the present paper is to provide a unified geometric framework for coordinated motion on Lie groups, proceeding as follows. (i) Coordination on Lie groups is defined from first principles of symmetry, distinguishing three variants: \emph{left-invariant}, \emph{right-invariant} and \emph{total} coordination. (ii) Expressing the conditions for coordination in the associated Lie algebra, a direct link is drawn between coordination on Lie groups and consensus in vector spaces. (iii) It is investigated how total coordination restricts compatible relative positions through a geometrically meaningful relation. These properties are independent of the system's dynamics. Going over to control laws, simplified first-order dynamics are assumed for individual agents, but underactuation is explicitly modeled; communication among agents is restricted to a reduced set of links that can possibly be directed and time-varying. (iv) Control laws based on standard vector space consensus algorithms are given that achieve the easier tasks of right-invariant coordination and fully actuated left-invariant coordination on general Lie groups, for any initial condition. (v) A general method is proposed to design control laws that achieve total coordination of fully actuated agents when the communication links among agents are undirected and fixed; extension to more general communication settings can be made along the lines of \cite{TAC2}. Total coordination design for fully actuated agents is a rather academic problem, but (vi) the proposed design methodology is then shown to apply to the practically most relevant problem of left-invariant coordination of underactuated agents. The proposed controller architecture consists of two steps, obtained by adding to the consensus algorithm a position controller derived from geometrical Lyapunov functions. The position controllers are directly linked to the double bracket flows of \cite{Brockett2} for gradient systems on adjoint orbits.

The power of the geometry is illustrated on $SO(3)$, $SE(2)$ and $SE(3)$ by analyzing the meaning of the geometric conditions for coordination, and by designing corresponding control laws with the proposed general methodology. The latter leads to controllers that have been previously proposed in the literature, but were derived based on intuitive arguments for these particular applications. In that sense, the novelty of the present paper is not in the expression of the coordinating control laws but in showing that they can be derived in a unifying and algorithmic manner with the proper geometric setting.\\

\paragraph*{Table of contents} The paper is organized as follows. Section 2 examines the geometric properties of coordination on Lie groups (contributions (i), (ii) and (iii)). Section 3 presents the control setting and basic control laws for right-invariant coordination and fully actuated left-invariant coordination (contribution (iv)). Sections 4 and 5 present control law design methods respectively for total coordination (contribution (v)) and for underactuated left-invariant coordination (contribution (vi)). Examples are treated at the end of Sections 2, 4 and 5.



\section{The geometry of coordination}

This section proposes definitions for coordination on Lie groups by starting from basic symmetry principles. It establishes conditions on the velocities for coordination and examines their implications. Except that the symmetries must be compatible, these developments are independent of the dynamics considered for the control problem. Notations are adapted from \cite{ArnoldMech}.


\subsection{Relative positions and coordination} Consider a swarm of $N$ ``agents'' evolving on a Lie group $G$, with $g_k(t)\in G$ denoting the position of agent $k$ at time $t$. Let $g_k^{-1}$ denote the group inverse of $g_k$, $L_h : g \mapsto hg$ denote left multiplication, and $R_h : g \mapsto gh$ right multiplication on $G$.\vspace{2mm}

\noindent \emph{Definition 1:} The \emph{left-invariant relative position} of agent $j$ with respect to agent $k$ is $\lambda_{jk} = g_k^{-1} g_j$. The \emph{right-invariant relative position} of agent $j$ with respect to agent $k$ is $\rho_{jk} = g_j \, g_k^{-1}$.\vspace{2mm}

Indeed, $\lambda_{jk}$ (resp. $\rho_{jk}$) is invariant under left (resp. right) multiplication: $(h g_k)^{-1}(h g_j) = g_k^{-1} g_j$ $\forall h \in G$. The left-/right-invariant relative positions are the \emph{joint invariants} associated to the left-/right-invariant action of $G$ on $G \times G ... \times G$ ($N$ copies).

The two different definitions of relative position lead to two different types of coordination; a third type is defined by combining them.
\vspace{2mm}

\noindent \emph{Definition 2:} \emph{Left-invariant coordination} (LIC) means constant left-invariant relative positions $\lambda_{jk}(t) = g_k^{-1}g_j$ --- resp. \emph{right-invariant coordination} (RIC) means constant right-invariant relative positions $\rho_{jk} = g_j g_k^{-1}$ --- for all pairs of agents $j,k$ in the swarm. \emph{Total coordination} (TC) means simultaneous left-invariant and right-invariant coordination: $g_k^{-1} g_j$ and $g_j g_k^{-1}$ are constant for all pairs of agents $j,k$ in the swarm.\vspace{2mm}

The present paper thus associates \emph{coordination} to fixed relative positions. In contrast, \emph{synchronization} is the situation where all agents are at the same point on $G$: $g_k = g_j$ $\forall j,k$; this is a very particular case of total coordination.


\subsection{Velocities and coordination}

Denote by $\al$ the Lie algebra of $G$, i.e. its tangent plane at the identity $e$. Denote by $[\; , \; ]$ the Lie bracket on $\al$. Let $L_{h*} : TG_g\rightarrow TG_{hg}$ and $R_{h*} : TG_g\rightarrow TG_{gh}$ for all $g\in G$ be the induced maps on tangent spaces corresponding to left- and right-multiplication $L_h$ and $R_h$ respectively. Let $Ad_g : \al \rightarrow \al$, $Ad_g = R_{g^{-1}*}L_{g*}$. \vspace{2mm}

\noindent \emph{Definition 4:} The left-invariant velocity $\xi^l_k \in \al$ and the right-invariant velocity $\xi^r_k \in \al$ of agent $k$ are defined by $\xi^l_k(\tau)=L_{g^{-1}(\tau)*}(\tfrac{d}{dt}g_k(t)\vert_{t=\tau})$ and $\xi^r_k(\tau)=R_{g^{-1}(\tau)*}(\tfrac{d}{dt}g_k(t)\vert_{t=\tau})$ respectively.\vspace{2mm}

The left-invariant (resp. right-invariant) velocity is such that $g_k(t)$ and $L_h g_k(t)$ (resp. $R_h g_k(t)$) have the same velocity $\xi^l_k(t)$ (resp. $\xi^r_k(t)$), for any fixed $h \in G$. Note the important equality 
\begin{equation}\label{xiL_xiR}
\xi^r_k=Ad_{g_k}\xi^l_k \; .
\end{equation}

\noindent \emph{Proposition 1:} Left-invariant coordination corresponds to equal right-invariant velocities $\xi^r_j = \xi^r_k$ $\forall j,k$. Right-invariant coordination corresponds to equal left-invariant velocities $\xi^l_j = \xi^l_k$ $\forall j,k$.\vspace{2mm}

\noindent \underline{Proof:} For $\lambda_{jk}$, $\tfrac{d}{dt} (g_k^{-1}g_j) = L_{g_k^{-1}*} \tfrac{d}{dt}g_j + R_{g_j*} \tfrac{d}{dt}g_k^{-1}$. But if $\tfrac{d}{dt}g_k = L_{g_k*} \xi_k^l$, then $\tfrac{d}{dt}g_k^{-1} = - L_{g_k^{-1}*} Ad_{g_k} \xi_k^l$. Thus $\tfrac{d}{dt}(g_k^{-1} g_j) = L_{g_k^{-1} g_j*} \xi^l_j - L_{g_k^{-1}*}R_{g_j*} Ad_{g_k} \xi_k^l = L_{g_k^{-1} g_j*} Ad_{g_j}^{-1} (Ad_{g_j} \xi^l_j - Ad_{g_k} \xi^l_k)$.
Since $L_{g_k^{-1} g_j*}$ and $Ad_{g_j^{-1}}$ are invertible, $\tfrac{d}{dt}(\lambda_{jk}) = 0$ is equivalent to $Ad_{g_j} \, \xi^l_j = Ad_{g_k} \, \xi^l_k$ or equivalently $\xi^r_j = \xi^r_k$. The proof for right-invariant coordination is strictly analogous. \hfill $\vartriangle$\\

Proposition 1 shows that coordination on the Lie group $G$ is equivalent to consensus in the vector space $\al$. The latter is a well-studied subject \cite{TsitsiklisThesis,MOREAU,MOREAU2,Tsitsiklis3,hendrickx1,olfati,ConsensusReview}. Total coordination requires \emph{simultaneous} consensus on $\xi_k^l$ and $\xi_k^r$; but the latter are not independent, they are linked through (\ref{xiL_xiR}) which depends on the agents' positions. \vspace{2mm}

\noindent \emph{Proposition 2:} Total coordination on a Lie group $G$ is equivalent to the following condition in the Lie algebra $\al$:
$$\forall k=1...N, \quad \xi^l_k = \xi^l \in \bigcap_{i,j} \mathrm{ker}(Ad_{\lambda_{ij}}-Id) \quad \text{or equivalently} \quad \xi^r_k = \xi^r \in \bigcap_{i,j} \mathrm{ker}(Ad_{\rho_{ij}}-Id)$$

\noindent \underline{Proof:} \RIC requires $\xi^l_k = \xi^l_j$ $\forall j,k$; denote the common value of the $\xi^l_k$ by $\xi^l$. Then \LIC requires $Ad_{g_k} \, \xi^l = Ad_{g_j} \, \xi^l$ $\Leftrightarrow$ $\xi^l = Ad_{\lambda_{jk}} \, \xi^l$ $\forall j,k$. The proof with $\xi^r$ is similar.\hfill $\vartriangle$\\

Proposition 2 shows that total coordination puts no constraints on the relative positions when the group is Abelian, since $Ad_{\lambda_{ij}} = Id$ in this case. In contrast, on a general Lie group, total coordination with non-zero velocity can restrict the set of possible relative positions as follows.\vspace{2mm}

\noindent \emph{Proposition 3:} Let $CM_{\xi} := \lbrace g \in G : Ad_g \, \xi = \xi \rbrace$.
\newline a. For every $\xi \in \al$, $CM_{\xi}$ is a subgroup of $G$.
\newline b. The Lie algebra of $CM_{\xi}$ is the kernel of $ad_{\xi} = [\xi,\, ]$, i.e. $\mathfrak{cm}_{\xi} = \lbrace \eta \in \al : [\xi,\eta] = 0 \rbrace$.\vspace{2mm}

\noindent \underline{Proof:} a. $Ad_e \, \xi = \xi$ $\forall \xi$ since $Ad_e$ is the identity operator. $Ad_{g} \, \xi = \xi$ implies $Ad_{g^{-1}} \, \xi = \xi$ by simple inversion of the relation. Moreover, if $Ad_{g_1} \, \xi = \xi$ and $Ad_{g_2} \, \xi = \xi$, then $Ad_{g_1 g_2} \, \xi = Ad_{g_1} \, Ad_{g_2} \, \xi = Ad_{g_1} \, \xi = \xi$. Thus $CM_{\xi}$ satisfies all group axioms and must be a subgroup of $G$.
\newline b. Let $g(t) \in CM(\xi)$ with $g(\tau) = e$ and $\tfrac{d}{dt}g(t)\vert_{\tau} = \eta$. Then $\eta \in \mathfrak{cm}_{\xi} =$ the tangent space to $CM_{\xi}$ at $e$. For constant $\xi$, $Ad_g(t) \xi = \xi$ implies $\tfrac{d}{dt}(Ad_g(t)) \xi = 0$, with the basic Lie group property $\tfrac{d}{dt}(Ad_g(t))\vert_{\tau} = ad_{\eta}$. Therefore $[\eta,\, \xi] = 0$ is necessary. It is also sufficient since, for any $\eta$ such that $[\eta,\, \xi] = 0$, the group exponential curve $g(t) = \mathrm{exp}(\eta t)$ belongs to $CM_{\xi}$. \hfill $\vartriangle$\\

$CM_{\xi}$ and $\mathfrak{cm}_{\xi}$ are called the isotropy subgroup and isotropy Lie algebra of $\xi$; these are classical mathematical objects in group theory \cite{MarsdenBook2}. From Propositions 2 and 3, one method to obtain a totally coordinated motion on Lie group $G$ is to
(1) choose $\xi^l$ in the vector space $\al$ and
(2) position the agents such that $\lambda_{jk} \in CM_{\xi^l}$ for a set of pairs $j,k$ corresponding to the edges of a connected undirected graph. Then indeed, $\xi_k^l = \xi^l$ $\forall k$ ensures \RIC, and $\lambda_{jk} \in CM_{\xi^l}$ implies $Ad_{\lambda_{jk}} \xi_k^l = \xi_k^l = \xi_j^l$ such that $\xi_k^r = Ad_{g_k} \xi_k^l = Ad_{g_j} \xi_j^l = \xi_j^r$ and \LIC is achieved as well. The same can be done with $\xi^r$ and the $\rho_{jk}$. Note that a swarm at rest ($\xi_k^l = \xi_k^r = 0$ $\forall k$) is always totally coordinated.\vspace{2mm}

\noindent \emph{Remark 1:} In many applications involving coordinated motion, reaching a particular \emph{configuration}, i.e. specific values of the relative positions, is also relevant. \cite{MY4} defines specific configurations as extrema of a cost function. Imposing relative positions in the (intersection of) set(s) $CM_{\xi}$ for some $\xi$ can be another way to classify specific configurations; unlike \cite{MY4}, it works for non-compact Lie groups. For compact groups, there seems to be no connection between configurations characterized through $CM_{\xi}$ and those defined by \cite{MY4}.\vspace{2mm}

\noindent \emph{Remark 2:} It is also possible, conversely, to consider fixed relative positions $\lambda_{jk}$ and characterize the set of velocities $\xi$ compatible with total coordination. For non-Abelian groups and a sufficiently large number $N$ of agents, this set generically reduces to $\xi = 0$. 


\subsection{Examples}
\vspace{3mm}

The special orthogonal groups $SO(n)$ and special Euclidean groups $SE(n)$, $n \geq 2$, are well characterized; their basic properties can even be found in control textbooks like \cite{Jurdjevic}. Left-invariant coordination for the particular examples of $SE(2)$ and $SE(3)$ was already formulated in Lie group notation in \cite{JandK2D,JandK3D}.

\paragraph*{$SO(3)$} A point $g$ on $SO(3)$ is represented by a 3-dimensional rotation matrix $Q$.
\begin{itemize}
\item Group multiplication, inverse and identity are the corresponding matrix operations. 
\item The Lie algebra $\mathfrak{so}(3)$ is the set of skew-symmetric $3 \times 3$ matrices $[\omega]^{\wedge}$, operations $L_{Q*} \xi$ and $R_{Q*} \xi$ are represented by $Q [\omega]^{\wedge}$ and $[\omega]^{\wedge} Q$ respectively. The invertible mapping
$$\left( \begin{array}{ccc}
0 & -\omega_{(3)} & \omega_{(2)} \\ \omega_{(3)} & 0 & -\omega_{(1)} \\ -\omega_{(2)} & \omega_{(1)} & 0
\end{array} \right) \in \mathfrak{so}(3) \qquad \begin{array}{c}
\underrightarrow{\; \; \; [\cdot]^{\vee}_{\phantom{kkk}}} \\ \overleftarrow{\; \; \; [\cdot]^{\wedge^{\phantom{kk}}}}
\end{array} \qquad
\left( \begin{array}{c}
\omega_{(1)} \\ \omega_{(2)} \\ \omega_{(3)}
\end{array}
\right) \in \mathbb{R}^3
$$
identifies $\mathfrak{so}(3) \ni [\omega]^{\wedge}$ with $\mathbb{R}^3 \ni \omega$.
\item With this identification, $Ad_{Q} \omega = Q \omega$ and $[\omega_k, \, \omega_j] = [\omega_k]^{\wedge} \omega_j = \omega_k \times \omega_j$ (vector product).
\item In the standard interpretation of $Q$ as rigid body orientation, $\omega^l$ and $\omega^r$ are the angular velocities expressed in body frame and in inertial frame respectively.
\item \LIC (equal $\omega_k^r$), \RIC (equal $\omega_k^l$) and \TC have a clear mechanical interpretation in this case.
\item For \TC with $\omega \neq 0$, $\mathfrak{cm}_{\omega} = \lbrace \lambda \omega : \lambda \in \mathbb{R}\rbrace$ and $CM_{\omega} = \lbrace$rotations around axis $\omega \rbrace$.

The dimension of $\mathfrak{cm}_{\xi^l}$ ($\Leftrightarrow$ of $CM_{\xi^l}$) is $1$; the agents rotate with the same angular velocity $\omega_k^r$ in inertial space and have the same orientation up to a rotation around $\omega_k^r$.\\
\end{itemize}

\paragraph*{$SE(2)$} The special Euclidean group in the plane $SE(2)$ describes all planar rigid body motions (translations and rotations). An element of $SE(2)$ can be written $g=(r,\theta) \in \mathbb{R}^2 \times S^1$ where $r$ denotes position and $\theta$ orientation.
\begin{itemize}
\item Group multiplication $g_1 g_2 = (r_1 + Q_{\theta_1} r_2, \theta_1 + \theta_2)$ where $Q_{\theta}$ is the rotation of angle $\theta$. Identity $e = (0, 0)$ and inverse $g^{-1} = (-Q_{-\theta} r, -\theta)$.
\item Lie algebra $\mathfrak{se}(2) = \mathbb{R}^2\times\mathbb{R} \ni \xi = (v, \omega)$. Operations $L_{g*} (v, \omega) = (Q_{\theta} v, \omega)$ and $R_{g*} (v, \omega) = (v + \omega Q_{\pi/2}r, \omega)$.
\item $Ad_g \, (v, \omega) = (Q_{\theta} v - \omega Q_{\pi/2} r, \omega)$ and $[(v_1,\omega_1),(v_2,\omega_2)] = (\omega_1 Q_{\pi/2} v_2 - \omega_2 Q_{\pi/2} v_1,0)$.
\item In the interpretation of rigid body motion, $v^l$ is the body's linear velocity expressed in body frame, $\omega^l = \omega^r =: \omega$ is its rotation rate. However, for $\omega \neq 0$, $v^r$ is not the body's linear velocity expressed in inertial frame. Instead, $s = \tfrac{-Q_{\pi/2}}{\omega} v^r$ is the center of the circle drawn by the rigid body moving with $\xi^r = (v^r, \omega)$. In \cite{SPL2005}, the intuitive argument to achieve coordination is to synchronize circle centers $s_k$; this actually means synchronizing right-invariant velocities $v^r_k$, $k=1...N$ ($\neq$ linear velocities expressed in inertial frame).
\item In \RIC, the agents move with the same velocity expressed in body frame (Figure 1, $r$). In \LIC, they move like a single rigid body: relative orientations and relative positions on the plane do not change (Figure 1, $l_1$ and $l_2$).
\item In \TC, the swarm moves like a single rigid body \emph{and} each agent has the same velocity expressed in body frame. Propositions 2 and 3 characterize $\mathfrak{cm}_{\xi^l}$ by $[\xi^l, \eta] = 0$ $\Leftrightarrow$ $\omega^l v_{\eta}= \omega_{\eta} v^l$ and $CM_{\xi^l}$ by $Ad_g \, \xi^l = \xi^l$ $\Leftrightarrow$ $(Q_{\theta}-\mathrm{Id}) v^l = \omega^ l Q_{\pi/2} r$. This leads to 3 different cases:
\begin{itemize}
\item[(o)] $\omega^l = v^l = 0$ $\Rightarrow$ $\mathfrak{cm}_{\xi^l} = \mathfrak{se}(2)$ and $CM_{\xi^l} = SE(2)$.
\item[(i)] $\omega^l=0$, $v^l \neq 0$ $\Rightarrow$ $\mathfrak{cm}_{\xi^l} = \lbrace (v, 0) : v \in \mathbb{R}^2 \rbrace$ and $CM_{\xi^l} = \lbrace (r, \, 0) : r \in \mathbb{R}^2 \rbrace$.
\item[(ii)] $\omega^l \neq 0$, any $v^l$ $\Rightarrow$ $\mathfrak{cm}_{\xi^l} = \lbrace (\tfrac{\omega}{\omega^l} v^l,\omega) : \omega \in \mathbb{R} \rbrace$. Define $C$, the circle of radius $\tfrac{\Vert v^l \Vert_2}{\vert \omega^l \vert}$ containing the origin, tangent to $v^l$ at the origin and such that $v^l$ and $\omega^l$ imply rotation in the same direction. Then solving $Ad_g \xi = \xi$ for $g$ and making a few calculations shows that $CM_{\xi^l} = \lbrace (r, \, \theta) : r \in C$ and $Q_{\theta} v^l$ tangent to $C$ at $r \rbrace$. This is consistent with an intuitive analysis of possibilities for circular motion with unitary linear velocity and fixed relative positions and orientations in the plane.
\end{itemize}
The dimension of $\mathfrak{cm}_{\xi^l}$ ($\Leftrightarrow$ of $CM_{\xi^l}$) is (o) 3, (i) 2 or (ii) 1. In case (o), the configuration is arbitrary but at rest. In case (i), the agents have the same orientation and move on parallel straight lines (Figure 1, $t_1$). In case (ii), they move on the same circle and have the same orientation with respect to their local radius (Figure 1, $t_2$).
\end{itemize}

\begin{figure}[ht]
\begin{center}
\setlength{\unitlength}{1mm}
\begin{picture}(155,75)
\put(0,15){\includegraphics[width=55mm]{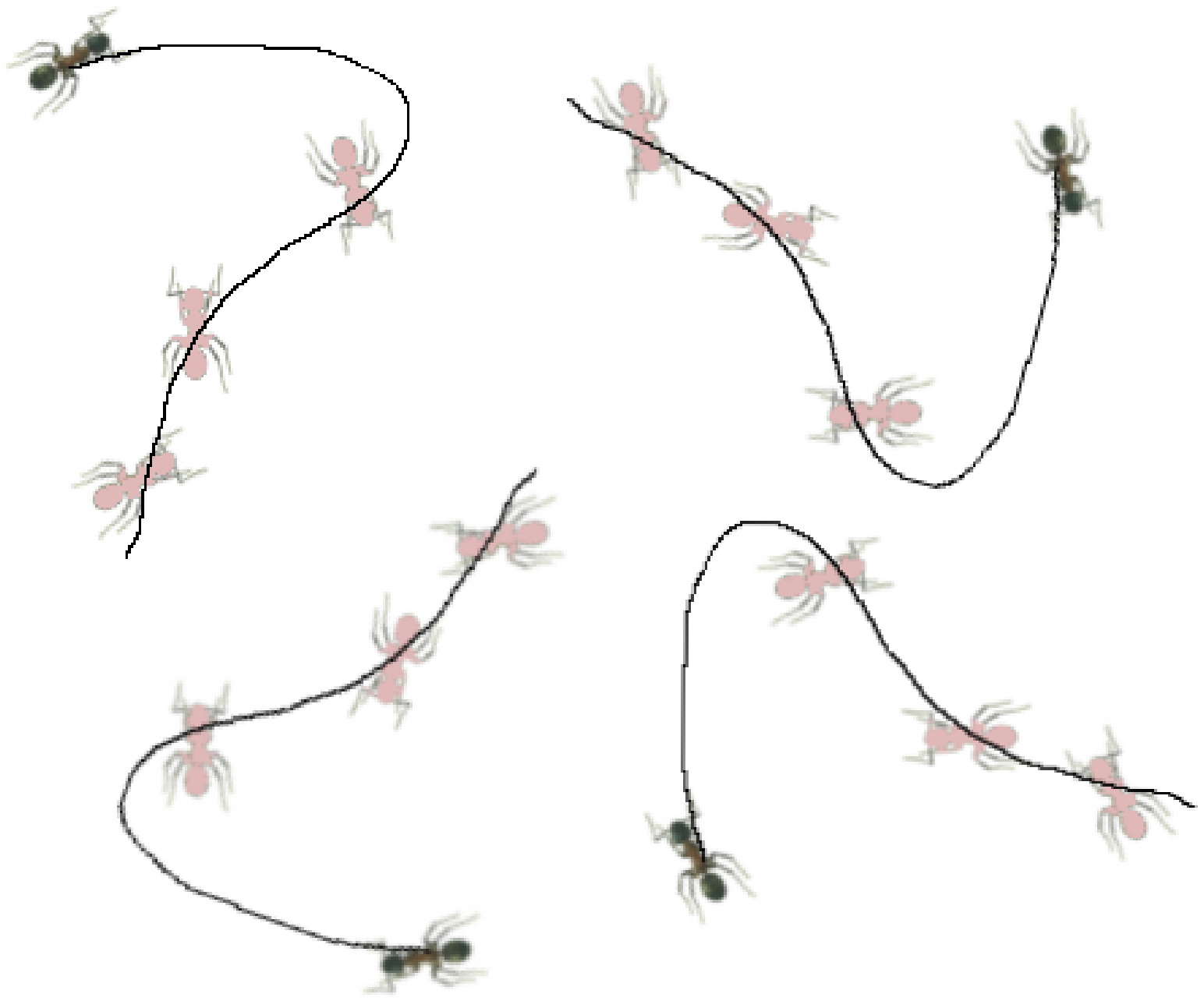}}
\put(0,15){$r$}
\put(65,45){\includegraphics[width=40mm]{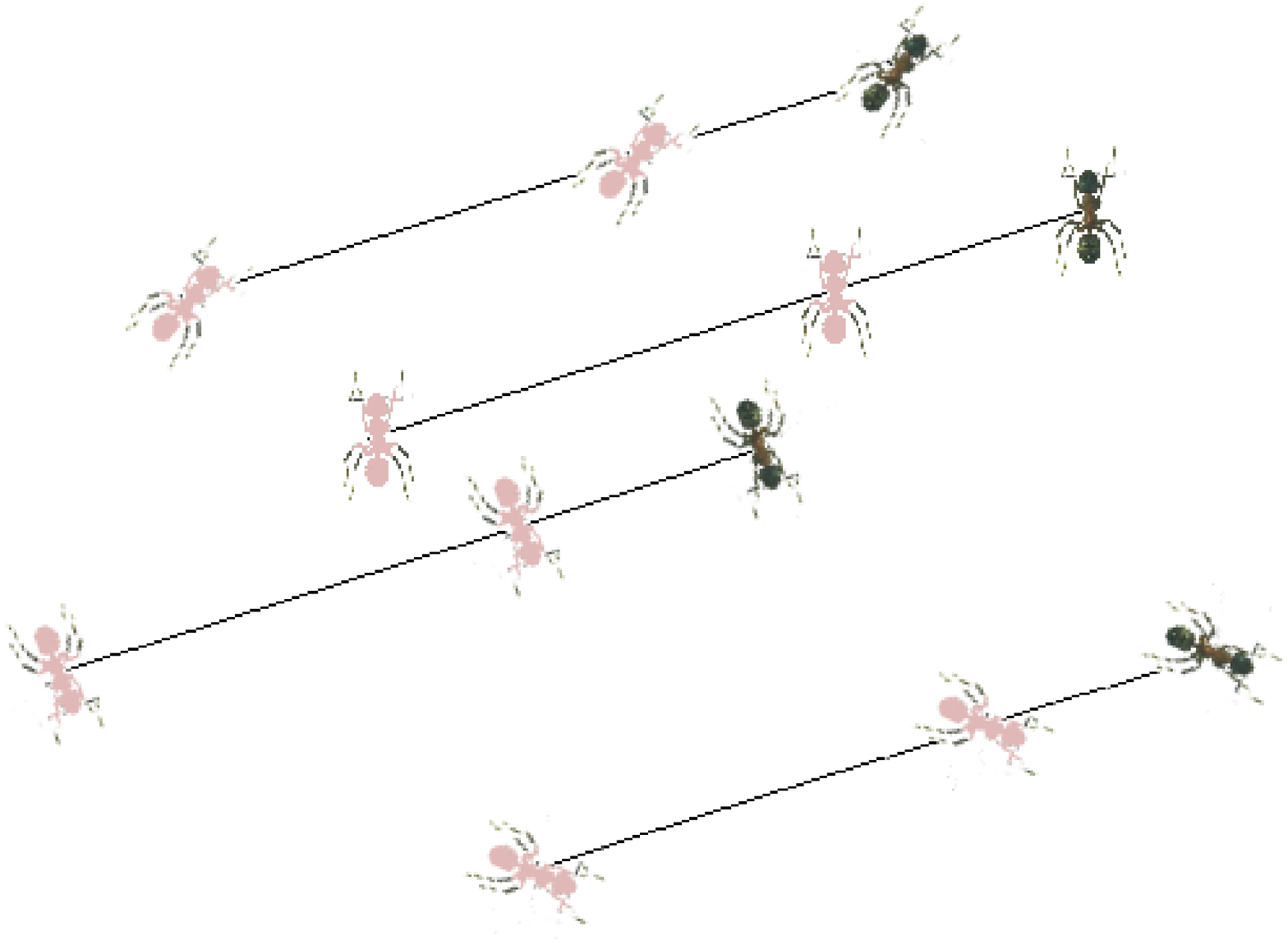}}
\put(65,45){$l_1$}
\put(65,0){\includegraphics[width=40mm]{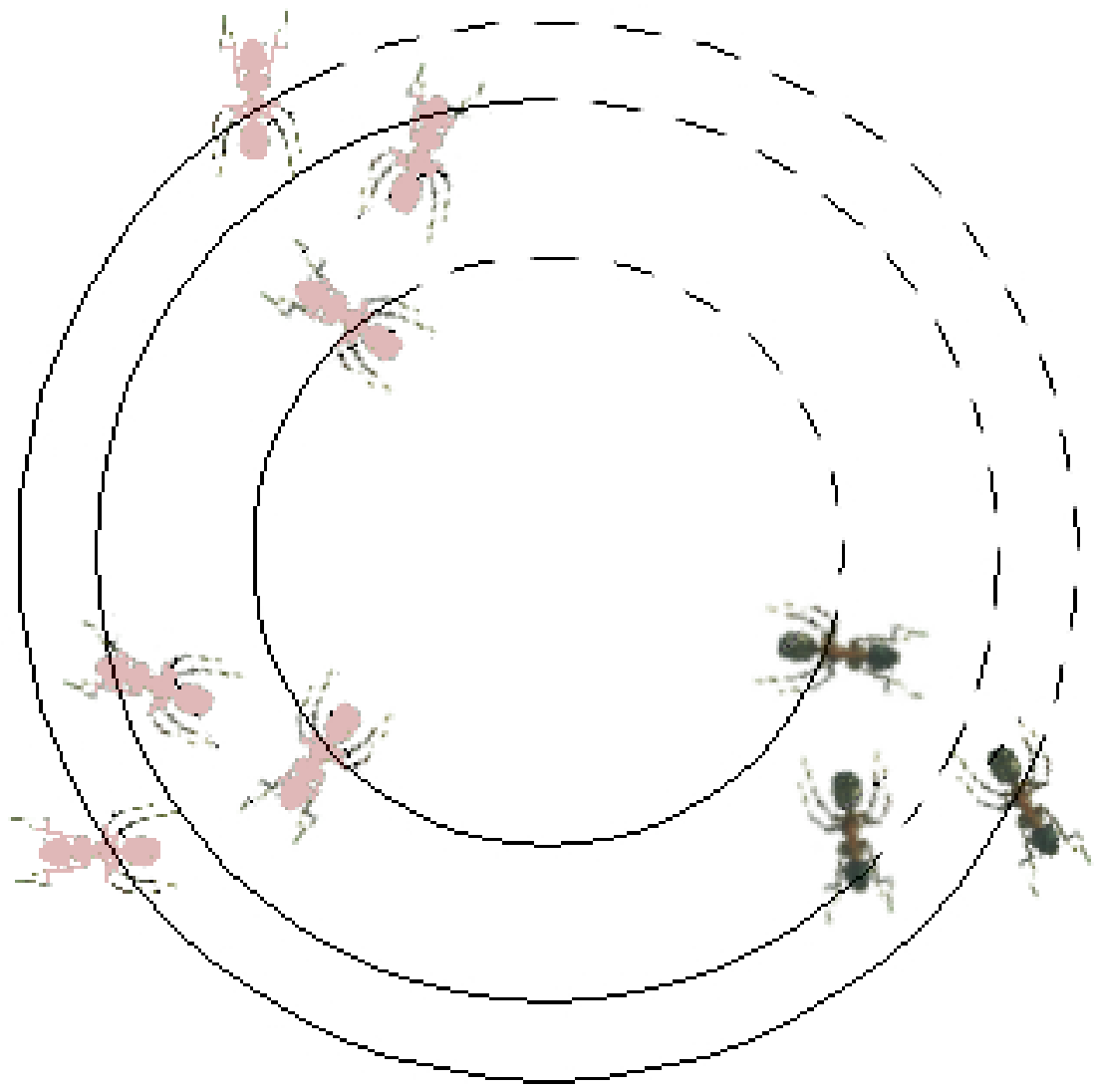}}
\put(65,0){$l_2$}
\put(115,45){\includegraphics[width=40mm]{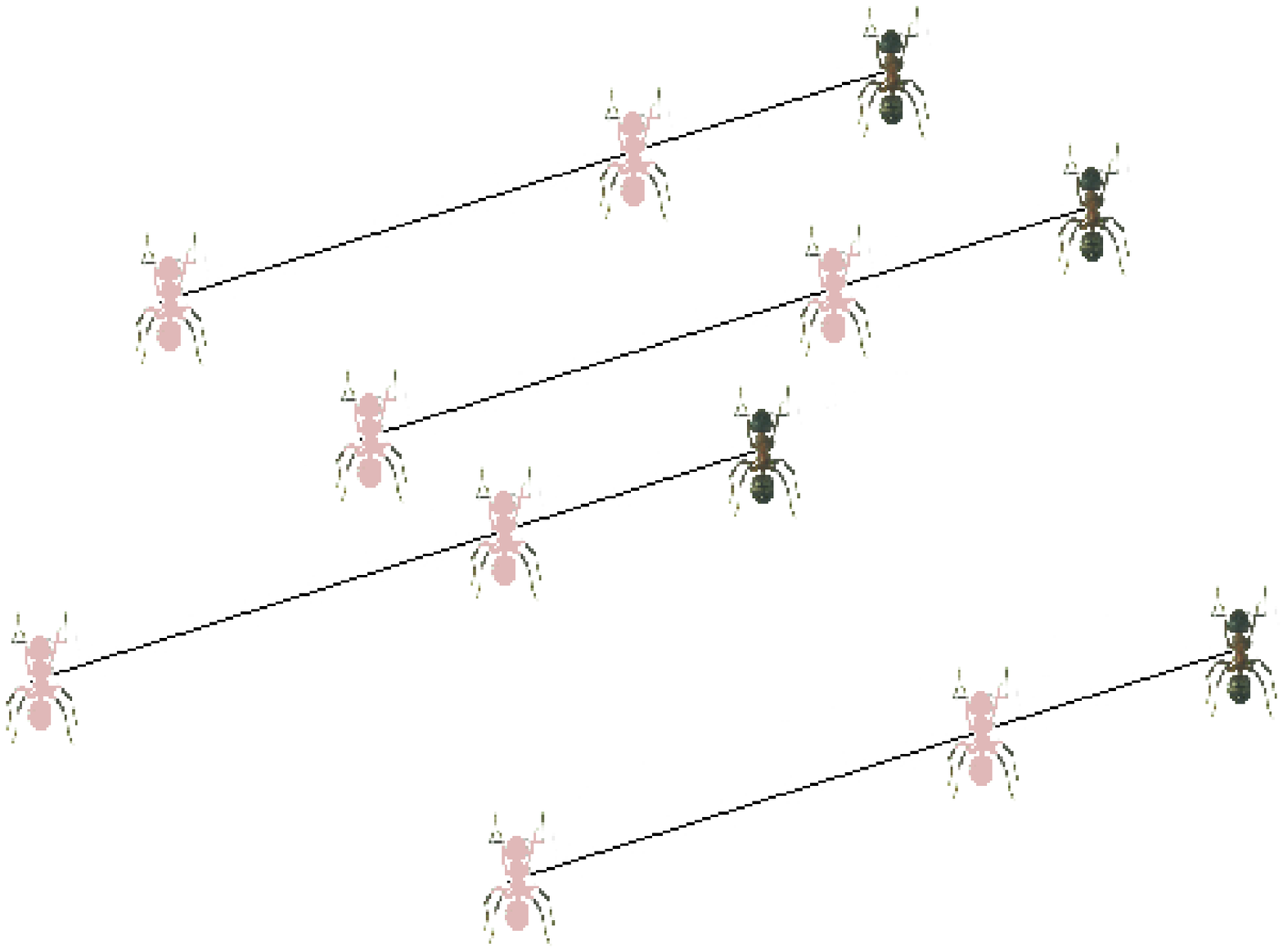}}
\put(118,45){$t_1$}
\put(115,0){\includegraphics[width=40mm]{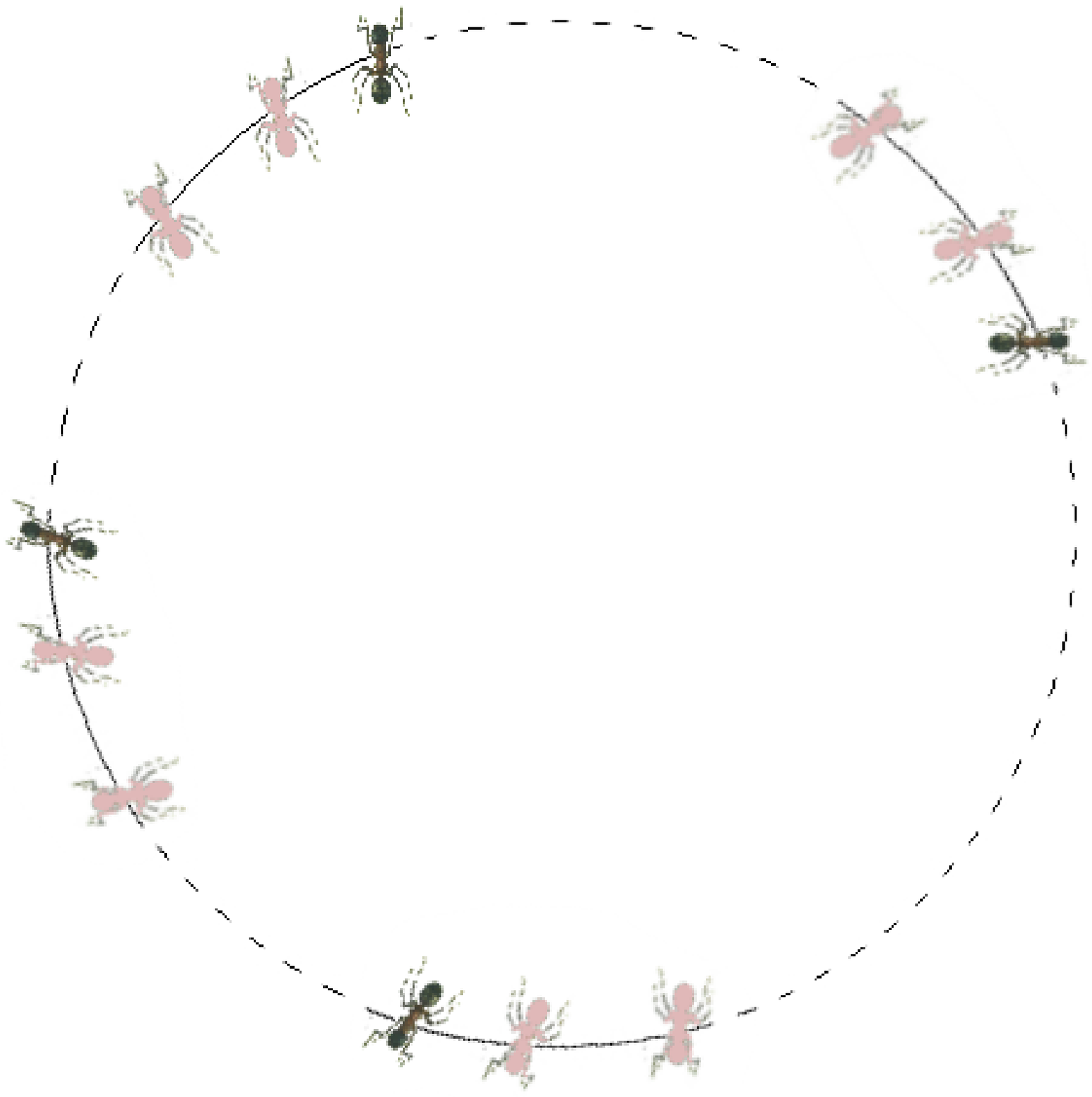}}
\put(118,0){$t_2$}
\end{picture}
\end{center}
\caption{Coordinated swarms (light color: intermediate positions and orientations in time). $r$: \RIC with varying velocity. $l_1$ and $l_2$: \LIC with $\omega_k = 0$ and $\omega_k \neq 0$ respectively. $t_1$ and $t_2$: \TC with $\omega_k = 0$ and $\omega_k \neq 0$ respectively.}
\end{figure}

\paragraph*{$SE(3)$} This group describes all $3$-dimensional rigid body motions (translations and rotations). An element of $SE(3)$ can be written $g=(r,Q) \in \mathbb{R}^3 \times SO(3)$, with $r$ denoting position and $Q$ orientation.
\begin{itemize}
\item $g_1 g_2 = (r_1 + Q_1 r_2, Q_1 Q_2)$, identity $e = (0, \mathrm{Id})$ and inverse $g^{-1} = (-Q^T r, Q^T)$.
\item Lie algebra $\mathfrak{se}(3) = \mathbb{R}^3\times\mathfrak{so}(3) \ni \xi = (v, [\omega]^{\wedge})$ is identified with $\mathbb{R}^3\times\mathbb{R}^3 \ni (v, \omega)$ with the same mapping as for $SO(3)$. Operations $L_{g*} (v, [\omega]^{\wedge}) = (Q v, Q [\omega]^{\wedge})$ and $R_{g*} (v, [\omega]^{\wedge}) = (\omega \times r + v, [\omega]^{\wedge} Q)$. As for $SO(3)$, symbol ``$\times$'' denotes vector product.
\item $Ad_g \, (v, \omega) = (Q v + r \times (Q \omega), Q \omega)$ and $[(v_1,\omega_1),(v_2,\omega_2)] = (\omega_1 \times v_2 - \omega_2 \times v_1, \omega_1 \times \omega_2)$.
\item In the interpretation of rigid body motion, left-invariant velocities $v^l$ and $\omega^l$ are the body's linear and angular velocity respectively, expressed in body frame; the right-invariant $\omega^r$ is the angular velocity expressed in inertial frame; for $\omega^l \neq 0$, there is no intuitive physical interpretation for the right-invariant $v^r$.
\item Similarly to $SE(2)$, the agents move in \RIC with the same velocity expressed in body frame and in \LIC with fixed relative orientations and relative positions, like a single rigid body.
\item In \TC, the swarm moves like a single rigid body \emph{and} each agent has the same velocity expressed in body frame. Propositions 2 and 3 lead to three different cases characterizing\linebreak
$\mathfrak{cm}_{\xi^l}$ which requires $[\xi^l, \eta] = 0$ $\Leftrightarrow \; \;$ $\omega^l \times \omega_{\eta} = 0$ and $\omega^l \times v_{\eta}=  \omega_{\eta} \times v^l$; \newline
$CM_{\xi^l}$ which requires $Ad_g \, \xi^l = \xi^l$ $\Leftrightarrow \; \;$ $Q \omega^l = \omega^l$ and $(Q-\mathrm{Id}) v^l = \omega^l \times r$.
\begin{itemize}
\item[(o)] $\omega^l = v^l = 0$ $\Rightarrow$ $\mathfrak{cm}_{\xi^l} = \mathfrak{se}(3)$ and $CM_{\xi^l} = SE(3)$.
\item[(i)] $\omega^l=0$, $v^l \neq 0$ $\Rightarrow$ $\mathfrak{cm}_{\xi^l} = \lbrace (\beta, \alpha v^l) : \beta \in \mathbb{R}^3,\; \alpha \in \mathbb{R} \rbrace$ and\newline $CM_{\xi^l} = \lbrace (r, \, Q) : r \in \mathbb{R}^3, \; Q$ characterizes rotation of axis $v^l \rbrace$.
\item[(ii)] $\omega^l \neq 0$, any $v^l$ $\Rightarrow$ $\mathfrak{cm}_{\xi^l} = \lbrace (\alpha v^l + \beta \omega^l, \alpha \omega^l) : \alpha, \, \beta \in \mathbb{R} \rbrace$ and $CM_{\xi^l} = \lbrace (r,\, Q) \in SE(3)$ describing left-invariant relative positions of agents that are on the same cylinder of axis $\omega^l$ and radius $\tfrac{\Vert v^l - (v^l) \cdot (\omega^l) / \Vert \omega^l \Vert \Vert}{\Vert \omega^l \Vert}$, with orientations differing around axis $\omega^l$ by an angle exactly equal to their relative angular position on the cylinder $\rbrace$. This is again obtained by solving for $g$ in $Ad_g \xi = \xi$ and making several basic computations; it is less obvious than for $SE(2)$ to find this result intuitively.
\end{itemize}
The dimension of $\mathfrak{cm}_{\xi^l}$ ($\Leftrightarrow$ of $CM_{\xi^l}$) is (o) 6, (i) 4 or (ii) 2. In case (o), the configuration is arbitrary but at rest. In case (i), the agents move on parallel straight lines and have the same orientation up to rotation around their linear velocity vector. In case (ii), for $v^l - (v^l) \cdot (\omega^l) / \Vert \omega^l \Vert \neq 0$, the agents draw helices of constant \emph{pitch} $\omega^l \cdot v^l = \omega^r \cdot v^r$ on the cylinder; in the special case $\omega^l \cdot v^l = 0$ the trajectories become circular (see figures in \cite{JandK3D,LUCA2}). In the degenerate situation $v^l - (v^l) \cdot (\omega^l) / \Vert \omega^l \Vert = 0$, all agents are on the rotation axis.
\end{itemize}



\section{Coordination as consensus in the Lie algebra}

\subsection{Control setting}

Left-invariant\footnote{A right-invariant system is equivalent, simply by redefining the group multiplications.} systems on Lie groups appear naturally in many physical systems, such as rigid bodies in space, and cart-like vehicles. Motivated by examples like 2-axes attitude control and steering control on $SE(2)$ or $SE(3)$, this paper considers a left-invariant dynamics with affine control of the type
\begin{equation}\label{ControlSetting}
\tfrac{d}{dt} g_k = L_{g_k*} \xi^l_k \qquad \; \text{with} \qquad \; \xi^l_k = a + B u_k \qquad , \; k=1...N \, ,
\end{equation}
where the Lie algebra $\al$ is identified with $\mathbb{R}^n$, $a \in \mathbb{R}^n$ is a constant  drift velocity, $B \in \mathbb{R}^{n \times m}$ has full column rank and specifies the range of the control term $u_k \in \mathbb{R}^m$. The set of all assignable $\xi_k^l$ is denoted $\C = \lbrace a + B u : u \in \mathbb{R}^m \rbrace$. Note that for fully actuated agents $m=n$, (\ref{ControlSetting}) boils down to $\tfrac{d}{dt} g_k = L_{g_k*} u_k$. Feedback control laws must be functions of variables which are compatible with the symmetries of the problem setting, i.e. left-invariant variables. In terms of left-invariant variables, \LIC corresponds to fixed (left-invariant) relative positions, while \RIC corresponds to equal (left-invariant) velocities.

In a realistic scalable setting, full communication between all agents cannot be assumed. The information flow among agents is modeled by restricting communication links among agents; $j \rightsquigarrow k$ denotes that $j$ sends information to $k$. The communication topology is associated to a graph $\G$. $\G$ is undirected if $k \rightsquigarrow j \, \Leftrightarrow \, j \rightsquigarrow k$. $\G$ is \emph{uniformly connected} (see \cite{MOREAU,MOREAU2}) if there exist an agent $k$ and durations $\delta > 0$ and $T > 0$ such that, $\forall t$, taking the union of the links appearing for at least $\delta$ in time span $[t,t+T]$, there is a directed path $k \rightsquigarrow a \rightsquigarrow b ... \rightsquigarrow j$ from $k$ to every other agent $j$.


\subsection{Right-invariant coordination}

Right-invariant coordination requires $\xi_k^l = \xi_j^l$ $\forall j,k$. In the setting (\ref{ControlSetting}), this simply implies to agree on equal $u_k$ $\forall k$; positions $\lambda_{jk}$ can evolve arbitrarily. This problem is solved by the classical vector space consensus algorithm \cite{TsitsiklisThesis,MOREAU2,Tsitsiklis3,hendrickx1,olfati,ConsensusReview}
\begin{equation}\label{L0}
\tfrac{d}{dt} \xi_k^l = {\textstyle \sum_{j \rightsquigarrow k}} \; (\xi_j^l - \xi_k^l) \qquad , \; k = 1...N \; ,
\end{equation}
that, using (\ref{ControlSetting}), translates into $\tfrac{d}{dt} u_k = \sum_{j \rightsquigarrow k} \; (u_j - u_k)$, and exponentially achieves $\xi_k^l = \xi_j^l$ $\forall j,k$ if $\G$ is uniformly connected. Agent $k$ relies on the left-invariant velocity $\xi_j^l$ of $j \rightsquigarrow k$; the initial values of $u_k$ can be chosen arbitrarily.\\

For a time-invariant and undirected communication graph $\G$, (\ref{L0}) is a gradient descent for the disagreement cost function $\; V_r = {\textstyle \sum_{k} \sum_{j \rightsquigarrow k}} \; \Vert \xi_k^l - \xi_j^l \Vert^2 \, , \; \;$ with the Euclidean metric in $\al$.


\subsection{Left-invariant coordination}\label{S:LIC}

Left-invariant coordination requires $\xi_k^r = \xi_j^r$ $\forall j,k$, which suggests to use
\begin{equation}\label{aux0}
\tfrac{d}{dt} \xi_k^r = {\textstyle \sum_{j \rightsquigarrow k}} \; (\xi_j^r - \xi_k^r) \qquad , \; k = 1...N \; .
\end{equation}
Using (\ref{xiL_xiR}) to rewrite (\ref{aux0}) in terms of the left-invariant variables yields
\begin{equation}\label{L1}
\tfrac{d}{dt} \xi^l_k = {\textstyle \sum_{j \rightsquigarrow k}} (Ad_{g_k^{-1} g_j} \, \xi^l_j - \xi^l_k) \qquad , \; k = 1...N
\end{equation}
thanks to $\tfrac{d}{dt}(Ad_{g_k} \xi^l_k) = Ad_{g_k} [\xi^l_k, \, \xi^l_k]=0$. To implement (\ref{aux0}), agent $k$ must know the relative position ${g_k^{-1} g_j}$ and velocity $\xi_j^l$ of $j \rightsquigarrow k$; the initial $u_k$ are still chosen arbitrarily.

The disagreement cost function $\; V_l = {\textstyle \sum_{k} \sum_{j \rightsquigarrow k}} \; \Vert Ad_{g_k} \xi^l_k - Ad_{g_j} \xi^l_j \Vert^2 \; \;$ associated to (\ref{aux0}) is not left-invariant in general (it involves positions $g_k$), so (\ref{L1}) cannot be a left-invariant gradient of $V_l$.

Nevertheless, let $\mathcal{G}_u$ be the subclass of compact groups with unitary adjoint representation, i.e. satisfying $\Vert \Ad_g \; \xi \Vert = \Vert \xi \Vert$ $\forall g \in G$ and $\forall \xi \in \al$ (for instance $SO(n)\in\mathcal{G}_u$). It is possible to define a bi-invariant (that is, left- and right-invariant) Riemannian metric on $G$ if and only if $G \in \mathcal{G}_u$ \cite{Pennec}. Using the Euclidean metric on left-invariant velocities, as in the present paper, comes down to using a left-invariant metric, in accordance with the left-invariant setting. If $G \in \mathcal{G}_u$, then this metric is bi-invariant, $ \; V_l = {\textstyle \sum_{k} \sum_{j \rightsquigarrow k}} \; \Vert \xi^l_k - Ad_{g_k^{-1} g_j} \xi^l_j \Vert^2 \; \;$ and for fixed undirected $\G$, (\ref{L1}) is a gradient descent for $V_l$.

A priori, the consensus algorithm (\ref{L1}) converges as (\ref{L0}). However, in contrast to (\ref{L0}), nothing guarantees that (\ref{L1}) can be implemented in an underactuated setting. At equilibrium, (\ref{L1}) requires 
\begin{equation}\label{compatibilityL:eq}
Ad_{\lambda_{jk}}  (a + B u_j) \, = \, a + B u_k \quad \forall j,k \; ,
\end{equation}
which may or may not hold depending on the relative positions of the agents. This issue motivates the further study of underactuated \LIC in Section 5. Similarly, total coordination requires simultaneous consensus on left- and right-invariant velocities. At equilibrium, this means that (\ref{compatibilityL:eq}) must hold with equal controls $u_k$, i.e. 
\begin{equation}\label{compatibilityT:eq}
Ad_{\lambda_{jk}}  (a + B u_k) \, = \, a + B u_k \quad \forall j,k \; ,
\end{equation}
which also puts constraints on the relative positions of the agents. For this reason, total coordination is further studied in Section 4.\\

In the following, it is assumed that the agents are controllable. Obviously, controllability is sufficient for coordination as it allows the agents to reach any position from any initial condition. However, it is not always necessary, as long as positions compatible with (\ref{compatibilityL:eq}) or (\ref{compatibilityT:eq}) are globally reachable; in particular, for Abelian groups $Ad_g = Id$ $\forall g$ so any positions satisfy (\ref{compatibilityL:eq}) and (\ref{compatibilityT:eq}); in that case, (underactuated) \LIC and \TC become trivial.



\section{Control design: fully actuated total coordination}


\subsection{Total coordination on general Lie groups}

Total coordination requires to satisfy two objectives, \LIC and \RIC, simultaneously. In a first step, assume that the agents have at their disposal a reference right-invariant velocity $\xi^r$ which they can track, such that \LIC is ensured if $\xi_k^l = Ad_{g_k}^{-1} \xi^r$ $\forall k$. It remains to simultaneously achieve \RIC, which, as previously shown, involves controlling relative positions. Writing
\begin{equation}\label{4:control_form}
\xi_k^l = \eta_k^l + \, q_k \; , \quad k=1...N \; ,
\end{equation}
where $\eta^l_k := Ad_{g_k}^{-1} \xi^r$, the question is how to design $q_k$ in order to achieve \TC. For fixed undirected communication graph $\G$, inspired by the cost function for \RIC, define
$$V_{tr}(g_1,g_2...g_N) = \tfrac{1}{2} {\textstyle \sum_k \sum_{j \rightsquigarrow k}} \; \Vert \eta_k^l - \eta_j^l\Vert^2$$
where $\Vert \; \Vert$ denotes Euclidean norm. $V_{tr}$ characterizes the distance from \RIC \emph{assuming that every agent implements $\xi_k^l = Ad_{g_k}^{-1} \xi^r$}. The time variation of $V_{tr}$ due to motion of $g_k$ is
\begin{equation}\label{DVt}
\tfrac{d}{dt} V_{tr} = 2 \; \, {\textstyle \sum_{k} \sum_{j \rightsquigarrow k}} \; (\eta^l_k - \eta^l_j) \cdot
[\eta_k^l,\, \xi_k^l]
\end{equation}
where $\cdot$ denotes the canonical scalar product in $\al$, because $\tfrac{d}{dt}(Ad_{g_k}^{-1}) \eta = -[\xi_k^l,\, Ad_{g_k}^{-1} \eta]$ $\forall \eta \in \al$. Thus if $q_k = 0$ then $\tfrac{d}{dt}V_{tr} = 0$, and a proper choice of $q_k$ should allow to decrease $V_{tr}$. Define\footnote{In fact, $\langle \, , \, \rangle$ expresses the effect of the Lie bracket on the dual space of $\al$, and is directly related to the \emph{coadjoint} representation of $G$; note however that in general, $\langle \, , \, \rangle$ does not satisfy the Lie bracket properties.} $\langle \, , \, \rangle$ such that $\xi_1 \cdot \langle \xi_2, \xi_3 \rangle + [\xi_1, \xi_2] \cdot \xi_3 = 0$ $\forall \xi_1,\xi_2,\xi_3 \in \al$; then (\ref{DVt}) rewrites $\tfrac{d}{dt} V_{tr} = 2 \; {\textstyle \sum_{k} \sum_{j \rightsquigarrow k}} \; \langle \eta_k^l, \eta^l_k - \eta^l_j \rangle \cdot q_k$ and the choice 
\begin{equation}\label{v_kfortotal}
q_k = -\langle \eta_k^l, \; {\textstyle \sum_{j \rightsquigarrow k}} \; (\eta^l_k - \eta^l_j) \rangle
\end{equation}
ensures that $V_{tr}$ is non-increasing along the solutions:
$$\tfrac{d}{dt}V_{tr} = -2 \; \, {\textstyle \sum_{k} \sum_{j \rightsquigarrow k}} \; \langle \eta_k^l, \; {\textstyle \sum_{j \rightsquigarrow k}} \; (\eta^l_k - \eta^l_j) \rangle^2 \; \leq 0 \, .$$

To obtain an autonomous, left-invariant algorithm for total coordination, it remains to replace the reference velocity $\xi^r$ by agent-related estimates $\eta_k^r$ on which the agents progressively agree. As the goal is to define a common right-invariant velocity in $\al$, it is natural to proceed as in Section \ref{S:LIC} and use the consensus algorithm 
\begin{equation}\label{4A:consensusr}
\tfrac{d}{dt}\eta_k^r = {\textstyle\sum_{j \rightsquigarrow k}} \, (\eta_j^r - \eta_k^r)
\end{equation}
which in terms of left-invariant velocities rewrites
\begin{equation}\label{4A:consensus}
\tfrac{d}{dt} \eta_k^l = {\textstyle \sum_{j\rightsquigarrow k}} (Ad_{\lambda_{jk}} \, \eta_j^l - \eta_k^l) \; - [\xi_k^l, \, \eta_k^l] \; , \quad k = 1...N \, .
\end{equation}

Thus the overall controller is the cascade of a consensus algorithm to agree on a desired velocity for \LIC, and a position controller designed to decrease a natural distance to \RIC. To implement the controller, agent $k$ must receive from communicating agents $j \rightsquigarrow k$ their relative positions $\lambda_{jk}$ and the values of their left-invariant \emph{auxiliary variables} $\eta_j^l$.

\begin{figure}[h!]
\begin{center}\setlength{\unitlength}{1mm}
\begin{picture}(110,15)
\put(0,0){\framebox(45,15){}}
\put(2,10){\LIC: agree on $\xi^r$}
\put(2,3){{\footnotesize vector space consensus in $\al$}}
\put(45,7.5){\vector(1,0){20}}
\put(65,0){\framebox(45,15){}}
\put(67,10){\RIC: agree on $Ad_{g_k}^{-1} \xi^r$}
\put(67,3){{\footnotesize Lyapunov-based control of $g_k$}}
\end{picture}
\end{center}
\caption{Total coordination as consensus on right-invariant velocity and Lyapunov-based control to right-invariant coordination.}
\end{figure}
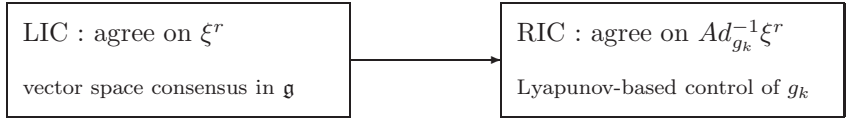

The following result characterizes the convergence properties of controller (\ref{4:control_form}),(\ref{v_kfortotal}),(\ref{4A:consensus}).\\

\noindent \emph{Theorem 1:} Consider $N$ fully actuated agents communicating on a fixed, undirected graph $\G$ and evolving on Lie group $G$ according to $\tfrac{d}{dt} g_k = L_{g_k*} \xi^l_k$ with controller (\ref{4:control_form}),(\ref{v_kfortotal}),(\ref{4A:consensus}).\begin{itemize}
\item[(i)] For any initial conditions $\eta_k^l(0)$, the $\eta_k^r(t)$ exponentially converge to $\overline{\eta^r} := \frac{1}{N} \sum_{k} \eta_k^r(0)$.
\item[(ii)] Define $\overline{V_{tr}}(g_1,g_2,...g_N) := \tfrac{1}{2} {\textstyle \sum_k \sum_{j \rightsquigarrow k}} \; \Vert Ad_{g_k}^{-1} \overline{\eta^r} - Ad_{g_j}^{-1} \overline{\eta^r} \Vert^2$. All solutions converge to the critical set of $\overline{V_{tr}}$. In particular, left-invariant coordination is asymptotically achieved.
\item[(iii)] Total coordination is (at least locally) asymptotically stable.
\end{itemize}\vspace{2mm}

\noindent \underline{Proof:} Regarding convergence, (\ref{4A:consensus}) is strictly equivalent to (\ref{4A:consensusr}). Therefore, (i) simply restates a well-known convergence result for consensus algorithms in vector spaces on fixed undirected graphs \cite{ConsensusReview}.

Since the $\eta_k^r$ converge, (\ref{4:control_form}),(\ref{v_kfortotal}) is an asymptotically autonomous system; the autonomous limit system is obtained by replacing $\eta_k^l = Ad_{g_k}^{-1} \overline{\eta^r}$. From the derivation of $q_k$ in (\ref{v_kfortotal}), the limit system is a gradient descent system for $\overline{V_{tr}}(g_1,g_2,...g_N)$; the latter is smooth because the adjoint representation is smooth. According to \cite{ChainRecurrent}, the $\omega$-limit sets of an asymptotically autonomous system correspond to the chain recurrent sets of the limit system. Moreover, from \cite{ChRec2} the chain recurrent set of a smooth gradient system is equal to the set of its critical points. Therefore the $\omega$-limit set of (\ref{4:control_form}),(\ref{v_kfortotal}) is equal to the set of critical points of $\overline{V_{tr}}$, which proves (ii). Total coordination $\overline{V_{tr}} = 0$ is locally asymptotically stable as it is a local (and global) minimum of $\overline{V_{tr}}$, which proves (iii). \hfill $\vartriangle$\\

Extensions to varying and directed $\G$ can be made with additional auxiliary variables along the lines of \cite{LUCA1,MY2,TAC2,LUCA2}: the algorithms define (estimate) a desired $\xi^l$ \emph{and} a desired $\xi^r$, which must be on the same adjoint orbit; cost functions for individual agents are used to ensure that they asymptotically implement the desired velocities. Sometimes intuition may be required to express everything in a left-invariant setting. These algorithms mostly overcome the problem of local minima different from \TC, which makes them useful for fixed undirected $\G$ as well.


\subsection{Total coordination on Lie groups with a bi-invariant metric}

When $G \in \mathcal{G}_u$, i.e. $G$ has a bi-invariant metric, the cost function $\; V_l = {\textstyle \sum_{k} \sum_{j \rightsquigarrow k}} \; \Vert Ad_{g_k} \xi^l_k - Ad_{g_j} \xi^l_j \Vert^2 \; \;$ can be used for left-invariant control design.

A natural idea in this context would be to combine the cost functions for \LIC and \RIC, writing $V_t = V_l + V_r$, and derive a gradient descent for $V_t$ of the form $\; \tfrac{d}{dt} \xi_k^l = f(\xi_k^l, \, \lbrace \xi_j^l, \, g_k^{-1} g_j : j \rightsquigarrow k \rbrace)$. However, simulations of the resulting control law for $SO(n)$ seem to always converge to $\xi_k^l = 0$ $\forall k$. A possible explanation for this behavior is that this strategy focuses on velocities, such that positions of the agents are not explicitly controlled, while it was shown in Section 2 that \TC at non-zero velocity involves restrictions on compatible positions.\\

Nevertheless, the existence of a bi-invariant metric offers the possibility to switch the roles of \LIC and \RIC in the method of Subsection 4.1, using a consensus algorithm to define a common left-invariant velocity for \RIC, and a cost function to drive positions to \LIC.

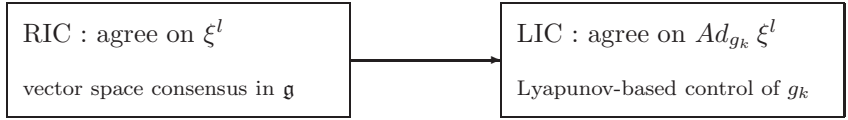
\begin{figure}[h!]
\begin{center}\setlength{\unitlength}{1mm}
\begin{picture}(110,15)
\put(0,0){\framebox(45,15){}}
\put(2,10){\RIC: agree on $\xi^l$}
\put(2,3){{\footnotesize vector space consensus in $\al$}}
\put(45,7.5){\vector(1,0){20}}
\put(65,0){\framebox(45,15){}}
\put(67,10){\LIC: agree on $Ad_{g_k}\, \xi^l$}
\put(67,3){{\footnotesize Lyapunov-based control of $g_k$}}
\end{picture}
\end{center}
\caption{Total coordination as consensus on left-invariant velocity and Lyapunov-based control to left-invariant coordination.}
\end{figure}

The \RIC consensus algorithm on auxiliary variables asymptotically define a common velocity $\xi^l$ by
\begin{equation}\label{eq:4B:cons}
\tfrac{d}{dt}\eta_k^l = {\textstyle\sum_{j \rightsquigarrow k}} \, (\eta_j^l - \eta_k^l) \; , \quad k=1...N \, .
\end{equation}
Then defining the cost function 
$$V_{tl}(g_1,g_2...g_N) = \tfrac{1}{2} {\textstyle \sum_k \sum_{j \rightsquigarrow k}} \; \Vert Ad_{g_k} \, \eta_k^l - Ad_{g_j} \, \eta_j^l\Vert^2 = \tfrac{1}{2} {\textstyle \sum_k \sum_{j \rightsquigarrow k}} \; \Vert  \eta_k^l - Ad_{g_k^{-1} g_j} \, \eta_j^l\Vert^2$$
for \LIC and proceeding as in the previous subsection, one obtains controller (\ref{4:control_form}) with
\begin{equation}\label{eq:4B:vk}
q_k = \langle \eta_k^l, \;{\textstyle \sum_{j \rightsquigarrow k}} \; (\eta^l_k - Ad_{g_k^{-1} g_j}\, \eta^l_j) \rangle \; .
\end{equation}\vspace{2mm}

\noindent \emph{Theorem 2:} Consider $N$ fully actuated agents communicating on a connected, fixed, undirected graph $\G$ and evolving on Lie group $G \in \mathcal{G}_u$ according to $\tfrac{d}{dt} g_k = L_{g_k*} \xi^l_k$ with controller (\ref{4:control_form}),(\ref{eq:4B:cons},(\ref{eq:4B:vk})).\begin{itemize}
\item[(i)] For any initial conditions $\eta_k^l(0)$, the $\eta_k^l(t)$ exponentially converge to $\overline{\eta^l} := \frac{1}{N} \sum_{k} \eta_k^l(0)$.
\item[(ii)] Define $\overline{V_{tl}}(g_1,g_2,...g_N) := \tfrac{1}{2} {\textstyle \sum_k \sum_{j \rightsquigarrow k}} \; \Vert Ad_{g_k} \overline{\eta^l} - Ad_{g_j} \overline{\eta^l} \Vert^2$. All solutions converge to the critical set of $\overline{V_{tl}}$. In particular, right-invariant coordination is asymptotically achieved.
\item[(iii)] Total coordination is (at least locally) asymptotically stable.
\end{itemize}\vspace{2mm}

\noindent \underline{Proof:} The proof is omitted because it is similar to the one of Theorem 1. \hfill $\vartriangle$\\

An advantage of Theorem 2 over Theorem 1 is that the control design can be directly extended to underactuated agents. Indeed, (\ref{eq:4B:cons}) defines a valid consensus velocity $\xi^l \in \C = \langle a + B u : u \in \mathbb{R}^m \rangle$ for underactuated agents provided that $\eta_k^l(0) \in \C$ $\forall k$. The only change is that $q_k$, instead of the exact gradient descent in (\ref{eq:4B:vk}), is its projection onto the control range of $B$:
$$\xi_k^l = a + B u_k = \eta_k^l + B \, B^T q_k$$
assuming without loss of generality that the columns of $B$ are orthonormal vectors. When $\xi^l$ is asymptotically defined with (\ref{eq:4B:cons}), the convergence argument for asymptotically autonomous systems must be extended to projections of the gradient system (\ref{eq:4B:vk}); a general proof of this technical issue is lacking in the present paper. It is the only reason to restrict Theorem 2 to fully actuated agents.\\

Brockett \cite{Brockett2} has developed a general double-bracket form for gradient algorithms on adjoint orbits of compact semi-simple groups, using the bi-invariant Killing metric. The connection with the present paper is obvious: once the consensus algorithm has converged, the gradient control for agent positions involves a cost function on the adjoint orbit of the common velocity $\overline{\eta^l}$ or $\overline{\eta^r}$. One example in \cite{Brockett2} involves minimizing the distance towards a subset of $\al$; a similar objective will be pursued in Section 5 of the present paper (but with a different class of subsets). A main difference of \cite{Brockett2} is its focus on the evolution of variables in $\al$, making abstraction of the underlying group, while in the present paper one actually controls the positions of (possibly underactuated) agents on $G$. If $G$ is a compact group and the bi-invariant Killing metric coincides with the left-invariant metric of the present paper, then $\langle \, , \, \rangle = -[ \, , \, ]$ and control (\ref{v_kfortotal}) for $g_k$ with $\eta_k^r = \xi^r$ fixed implies that $\eta_k^l$ follows the double bracket flow
\begin{equation}\label{DoubleBracket}
\tfrac{d}{dt} \eta_k^l = [\eta_k^l, \, [\eta_k^l,\, {\textstyle \sum_{j \rightsquigarrow k}} (\eta_k^l - \eta_j^l)]] \; .
\end{equation}
This is the case among others for the following example in $SO(3)$.

\subsection{Example: Total coordination in $SO(3)$}

Control laws for coordination in $SO(3)$ abound in the literature --- see among others the papers about satellite attitude control mentioned in the Introduction. Total coordination on $SO(3)$ requires aligned rotation axes, and thus synchronizes satellite attitudes up to their phase around the rotation axis.

The compact group $SO(3)$ has a bi-invariant metric, so the approach of Section 4.2 can be applied. Algorithm (\ref{eq:4B:cons}) is used verbatim, with $\eta_k^l \in \mathbb{R}^3$ the auxiliary variable associated to angular velocity $\omega_k^l$. As mentioned before equation (\ref{DoubleBracket}), $\langle \, , \, \rangle = -[ \, , \, ]$ on $SO(3)$. Thus in the fully actuated case, (\ref{4:control_form}),(\ref{eq:4B:vk}) lead to
\begin{equation}\label{eq:4:SO3a}
\omega_k^l = \eta_k^l \; + \eta_k^l \times ({\textstyle \sum_{j \rightsquigarrow k}} \; Q_k^T Q_j \eta^l_j ) \; , \quad k = 1...N \; .
\end{equation}
Theorem 2 can be strengthened as follows for specific graphs.\vspace{2mm}

\noindent \emph{Proposition 4:} If $\G$ is a tree or complete graph, \TC is the only asymptotically stable limit set.
\vspace{2mm}

\noindent \underline{Proof:} According to Theorem 2, it remains to show that \TC is the only local minimum of $V_{tl}$. Fixing $\eta_k^l = \omega^l$ $\forall k$, critical points of $V_{tl}$ correspond to
\begin{equation}\label{SO3_cond}
(Q_k \omega^l) \times ({\textstyle \sum_{j \rightsquigarrow k}} Q_j \omega^l) = 0 \quad \forall k \, .
\end{equation}
For the tree, start with the leaves $c$. Then $(Q_c \omega^l) \times (Q_p \omega^l) = 0$ where $p$ is the parent of $c$. As a consequence, (\ref{SO3_cond}) for the parent becomes $(Q_p \omega^l) \times (Q_{pp}^T \omega^l) = 0$ where $pp$ is the parent of $p$. Using this argument up to the root, all $(Q_k \omega^l)$ must be parallel. If the agents are partitioned in two anti-aligned groups, then moving those groups towards each other decreases $V_{tl}$; thus $V_{tl} = 0$ is the only local minimum. For the complete graph, (\ref{SO3_cond}) becomes $(Q_k \omega^l) \times \psi = 0$ $\forall k$, where $\psi = \sum_j Q_j \omega^l$. This implies either that all $Q_k \omega^l$ must be parallel or that $\psi = 0$. In the first case, further discussion is as for the tree. Rewriting $V_{tl} = N^2 \; \Vert \omega^l \Vert^2 - \tfrac{1}{2} \psi \cdot \psi$ shows that $\psi = 0$ corresponds to a maximum of $V_t$. \hfill $\vartriangle$\\

It is straightforward to adapt (\ref{eq:4:SO3a}) for underactuated agents; a popular underactuation on $SO(3)$ is to consider 2 orthogonal axes of allowed rotations $\mathbf{e}_1$ and $\mathbf{e}_2$, either controlling both rotation rates, i.e. $\omega_k^l = u_1 \mathbf{e}_1 + u_2 \mathbf{e}_2$, or imposing a fixed rotation rate around one axis, i.e. $\omega_k^l = \mathbf{e}_1 + u_2 \mathbf{e}_2$. Both cases are controllable \cite{Jurdjevic}, so the Jurdjevic-Quinn theorem ensures local asymptotic stability of \TC, if $\eta_k^l = \overline{\eta^l}$ $\forall k$ is fixed in advance or agreed on in finite time. A formal convergence proof for the asymptotically autonomous case where the $\eta_k^l$ follow (\ref{eq:4B:cons}) is currently missing.



\section{Control design: underactuated left-invariant coordination}

Total coordination may appear as a rather academic example, whose motivation in applications is not clear. However, the methodology developed in Section 4 for \TC control design is instrumental to achieve left-invariant coordination of underactuated agents. The latter is well motivated by practical applications. Here the role of the cost function is no longer to add a second level of coordination, but to fulfill the underactuation constraints. Unlike the academic example of \TC, the present section explicitly considers \LIC control design in the most general setting of underactuated agents as well as possibly directed and time-varying interconnection graph $\G$.

\subsection{Left-invariant coordination of underactuated agents}

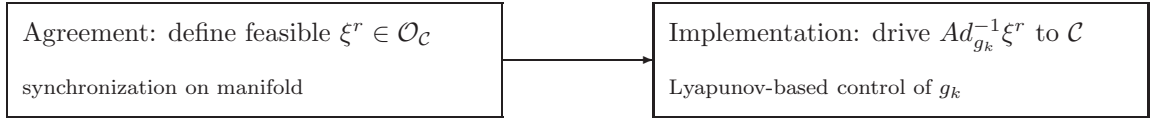
\begin{figure}[h!]
\begin{center}\setlength{\unitlength}{1mm}
\begin{picture}(150,15)
\put(0,0){\framebox(65,15){}}
\put(2,10){Agreement: define feasible $\xi^r \in \mathcal{O}_{\C}$}
\put(2,3){{\footnotesize synchronization on manifold}}
\put(65,7.5){\vector(1,0){20}}
\put(85,0){\framebox(65,15){}}
\put(87,10){Implementation: drive $Ad_{g_k}^{-1} \xi^r$ to $\C$}
\put(87,3){{\footnotesize Lyapunov-based control of $g_k$}}
\end{picture}
\end{center}
\caption{Underactuated left-invariant coordination as constrained consensus on right-invariant velocity and Lyapunov-based control to left-invariant coordination.}\label{MyFig}
\end{figure}

The control design of underactuated left-invariant coordination is decomposed in the two steps illustrated in Figure \ref{MyFig}. In a way analogous to the total coordination design of Section 4.1, a feasible right-invariant velocity is determined by a consensus algorithm. The corresponding left-invariant velocity is enforced by a Lyapunov-based feedback that decreases the distance of the consensus velocity to $\C = \lbrace a + B u : u \in \mathbb{R}^m \rbrace$.

The consensus algorithm must enforce a feasible right-invariant velocity, that is a vector $\xi^r$ in the set
$$O_{\C} := \lbrace Ad_g \xi : \xi \in \C \text{ and } g \in G \rbrace \, .$$
If $O_{\C}$ is convex, then it is sufficient to initialize the consensus algorithm (\ref{4A:consensus}) with $\eta_k^l(0) \in \C$. When $O_{\C}$ is not convex, the consensus algorithm must be adapted and the present paper has no general method. Strategies inspired from \cite{MY4} for compact homogeneous manifolds may be helpful, as illustrated in the example below.

Assuming a known feasible right-invariant velocity $\xi^r$, the design of a Lyapunov based control to left-invariant coordination proceeds similarly to Section 4.1.

Define $d(\eta,\, \C)$ to be the Euclidean distance in $\al$ from $\eta$ to the set $\C$.  Let $\Pi_{\C}(\eta)$ be the projection of $\eta$ on $\C$; since $\C$ is convex, $\forall \eta$~~$\Pi_{\C}(\eta)$ is the unique point in $\C$ such that $d(\eta,\, \C) = d(\eta,\, \Pi_{\C}(\eta)) =: \Vert \eta - \Pi_{\C}(\eta) \Vert$. Following the same steps as in Section 4.1, define $\eta_k^l := Ad_{g_k}^{-1} \xi^r$. Writing
\begin{equation}\label{eq:5:control:form}
\xi_k^l = a + B u_k = \Pi_{\C}(\eta_k^l) + B q_k  \; , \quad k=1...N \; ,
\end{equation}
the task is to design $q_k \in \mathbb{R}^m$ such that asymptotically, $g_k$ is driven to a point where $\eta_k^l \in \C$ and $q_k$ converges to $0$; this would asymptotically ensure \LIC. For each individual agent $k$, write the cost function
\begin{equation}\label{5:Vk}
V_k(g_k) = \tfrac{1}{2}\Vert Ad_{g_k}^{-1} \xi^r - \Pi_{\C}(Ad_{g_k}^{-1} \xi^r) \Vert^2=\tfrac{1}{2}\Vert\eta_k^l - \Pi_{\C}(\eta_k^l) \Vert^2
\end{equation}
where $\Vert \; \Vert$ denotes Euclidean norm. $V_k$ characterizes the distance of $\eta_k^l$ from $\C$, that is the distance from \LIC \emph{assuming that every agent implements $\xi_k^l = \Pi_{\C}(Ad_{g_k}^{-1} \xi^r)$}. The time variation of $V_k$ due to motion of $g_k$ is
\begin{equation}\label{DVkt}
\tfrac{d}{dt} V_k = (\eta_k^l - \Pi_{\C}(\eta_k^l)) \cdot [\eta_k^l,\, \Pi_{\C}(\eta_k^l) + B q_k]
\end{equation}
where $\cdot$ denotes the canonical scalar product in $\al$. Going further along the lines of Section 4.1 requires to assume that the control setting (pair $a$, $B$) and Lie algebra structure are such that $\forall \eta \in O_{\C}$, it holds $(\eta - \Pi_{\C}(\eta)) \cdot [\eta,\, \Pi_{\C}(\eta)] \leq 0$; then (\ref{DVkt}) implies $\tfrac{d}{dt} V_k \leq  f(\eta_k^l) \cdot q_k$ for some continuous $f:\al \rightarrow \mathbb{R}^m$ and a natural control is
\begin{equation}\label{v_kforleft}
q_k = -f(\eta_k^l) \; , \quad k=1...N \, .
\end{equation}
Note that when $O_{\xi^r} \subseteq \C$, the position control algorithm is unnecessary and vanishes, yielding simply
$\xi_k^l = Ad_{g_k}^{-1} \xi^r$ $\forall t$.

The overall controller is the cascade of a consensus algorithm to agree on a desired velocity for \LIC, and a position controller designed from a natural Lyapunov function to reach positions compatible with underactuation constraints and actually achieve \LIC. To implement the controller, agent $k$ must get from other agents $j \rightsquigarrow k$ their relative positions $\lambda_{jk}$ and the values of their left-invariant \emph{auxiliary variables} $\eta_j^l$. Since agents only interact through the consensus algorithm, not through the cost function, a connected fixed undirected graph is not required: $\G$ can be directed and time-varying, as long as it remains uniformly connected.

A general characterization of the behavior of solutions of the closed-loop system is more difficult here because the position controller is not a gradient anymore. The following result involves assumptions that can be readily checked for any particular case.\vspace{2mm}

\noindent \emph{Theorem 3:} Consider $N$ underactuated agents communicating on a uniformly connected graph $\G$ and evolving on Lie group $G$ according to $\tfrac{d}{dt} g_k = L_{g_k*} \xi^l_k$ with controller (\ref{eq:5:control:form}),(\ref{v_kforleft}), assuming that $\forall \eta \in O_{\C}$, it holds $(\eta - \Pi_{\C}(\eta)) \cdot [\eta,\, \Pi_{\C}(\eta)] \leq 0$. Assume that an appropriate consensus algorithm drives the arbitrarily initiated $\eta_k^l$, $k=1...N$, such that they exponentially agree on $Ad_{g_k} \eta_k^l \rightarrow \xi^r \in O_{\C}$ $\forall k$, independently of the agent motions $g_k(t)$.
\begin{itemize}
\item[(i)] If the agents are controllable, \LIC is locally asymptotically stable.
\item[(ii)] If, for any fixed $\eta_k^r = \xi^r$, bounded $V_k$ implies bounded $\eta_k^l$, and $f(\eta_k^l) \rightarrow 0$ implies $g_k \rightarrow \lbrace g : f(Ad_g^{-1} \eta_k^r) = 0 \rbrace$, then all agent trajectories on $G$ converge to the set where $f(Ad_{g_k}^{-1} \xi^r) = 0$.
\end{itemize}\vspace{2mm}

\noindent \underline{Proof:} The overall system is a cascade of the exponentially stable consensus algorithm and position controller (\ref{eq:5:control:form}),(\ref{v_kforleft}) which is decoupled for the individual agents. Assumptions $\tfrac{d}{dt} V_k \leq  f(\eta_k^l) \cdot q_k$ and (\ref{v_kforleft}) exactly mean that $V_k(g_k)$ is non-increasing along the closed-loop solutions. Therefore, if the agents are controllable, Jurdjevic-Quinn theorem \cite{JurdjevicQuinn} implies local asymptotic stability of the local minimum $V_k = 0$ $\forall k$ for the position controller. Then the overall system is the cascade of an exponentially stable system and a system for which $V_k = 0$ $\forall k$ is locally asymptotically stable. Standard arguments on cascade systems (see e.g. \cite{CascSontag,RodolpheControlBook}) allow to conclude that $V_k = 0$ $\forall k$ is locally asymptotically stable for the overall system; this proves (i).

To prove (ii), first consider the case where $\eta_k^r = \xi^r$ constant $\forall k$. Then $V_k$ can only decrease, and since it is bounded from below it tends to a limit; therefore $\tfrac{d}{dt}V_k$ is integrable in time for $t \rightarrow + \infty$. For the same reason, $V_k$ is bounded, so according to the assumption for (ii) $\eta_k^l$ is bounded as well; then $\tfrac{d^2}{dt^2}V_k$, which is a continuous function of $\eta_k^l$, is bounded as well for the closed-loop system, such that $\tfrac{d}{dt}V_k$ is uniformly continuous in time for $t \rightarrow + \infty$. Barbalat's Lemma implies that $\tfrac{d}{dt}V_k$ converges to $0$, which implies that $f(\eta_k^l)$ converges to 0, concluding the proof. Now in fact $\eta_k^r$ is not constant but exponentially converges to the constant value $\xi^r$ $\forall k$. But this changes nothing to the fact that $V_k$ tends to a finite limit and $\tfrac{d^2}{dt^2}V_k$ is bounded, so the same argument applies. \hfill $\vartriangle$\\

Condition $\tfrac{d}{dt} V_k \leq  f(\eta_k^l) \cdot q_k$ in Theorem 3 is not always true when $a \neq 0$; however, it is often satisfied in practice, as for steering control of rigid bodies in the following example. For this example, it is also possible to slightly improve Theorem 3 by showing that \LIC is the \emph{only stable} limit set.

\subsection{Example: Steering control on $SE(3)$}

Left-invariant coordination on $SE(3)$ under steering control is studied in \cite{LUCA2,LUCA3}. The present section shows how the algorithms of \cite{LUCA2} follow from the present general framework.

Using the notations of Section 2.3, the position and orientation of a rigid body in 3-dimensional space is written $(r_k,Q_k) =: g_k$, which is an element of the Special Euclidean group $SE(3)$; group multiplication is the usual composition law for translations and rotations, see Section 2.3. Then requiring agents to ``move in formation'', i.e. such that the relative position and heading of agent $j$ with respect to agent $k$ is fixed in the reference frame of agent $k$, $\forall j,k$, is equivalent to requiring left-invariant coordination. Moreover, since linear and angular velocity in body frame correspond to the components $(v^l_k,\omega^l_k)$ of $\xi^l_k$, the problem of controlling each agent \emph{in its own frame} with feedback involving \emph{relative} positions and orientations of other agents only, fits the left-invariant problem setting described in Section 3. The constraint of \emph{steering control} --- i.e. fixed linear velocity in agent frame $v^l_k = \mathbf{e}_1$ --- implies (\ref{ControlSetting}) of the form
$$\xi_k^l = a + B u_k = (\mathbf{e}_1, u_k) \quad \Rightarrow \quad \C = (\mathbf{e}_1, \mathbb{R}^3)\; .$$
Steering controlled agents on $SE(3)$ are controllable \cite{Jurdjevic}.

Following the method of Section 5.1, write auxiliary variables $\eta_k^l = (\eta_{v \, k}^l, \eta_{\omega \, k}^l)$; then $\Pi_{\C}(\eta_k^l) = (\mathbf{e}_1, \eta_{\omega \, k}^l)$, cost function $V_k = \tfrac{1}{2} \Vert \eta_{v \, k}^l - \mathbf{e}_1 \Vert^2$ and straightforward calculations show that (\ref{DVkt}) becomes $\tfrac{d}{dt} V_k = (\eta_{v \, k}^l \times \mathbf{e}_1) \cdot q_k$. This means that $(\eta - \Pi_{\C}(\eta)) \cdot [\eta,\, \Pi_{\C}(\eta)] = 0$ and $f_{\eta_k^l} = (\eta_{v \, k}^l \times \mathbf{e}_1)$. Then (\ref{eq:5:control:form}),(\ref{v_kforleft}) yield the controller
\begin{equation}\label{SE3:uk}
u_k = \eta_{\omega \, k}^l + \mathbf{e}_1 \times \eta_{v \, k}^l \; , \quad k=1...N \, .
\end{equation}
This is the same control law as derived in \cite{LUCA2} from intuitive arguments. If an appropriate consensus algorithm is built, then all assumptions of Theorem 3 hold, implying local asymptotic stability of 3-dimensional ``motion in formation'' with steering control (\ref{SE3:uk}); in fact, \cite{LUCA2} slightly improves Theorem 3 by also showing that \emph{globally}, \LIC is the only stable limit set.

It remains to design a consensus algorithm for the $\eta_k^l$. For this, two cases are distinguished, as in \cite{LUCA2}: linear motion $\omega^r = 0$ and helicoidal (of which a special case is circular) motion $\omega^r \neq 0$. The first case (almost) never appears from a consensus algorithm with arbitrary $\eta_k^l(0)$; it can however be imposed by $\eta_{\omega \, k}^l(0) = 0$ $\forall k$, which will then remain true $\forall t \geq 0$, in order to stabilize a coordinated motion in straight line.
\begin{itemize}
\item If $\eta_{\omega \, k}^l = 0$ (linear motion), then $\eta_{v \, k}^l = Q_k^T \eta_{v \, k}^r$ and $O_{(\mathbf{e}_1, 0)} = \lbrace (\lambda, 0) \in \mathbb{R}^3 \times \mathbb{R}^3 : \Vert \lambda \Vert = 1 \rbrace$. Agreement on $v^r$ in the unit sphere can be achieved following \cite{MY4}, just achieving consensus in $\mathbb{R}^3$ and normalizing; in fact normalizing is not even necessary, as it would just change the gain in (\ref{SE3:uk}). This leads to
\begin{equation}\label{SE3:cons1}
\tfrac{d}{dt}\eta_{v \, k}^l = {\textstyle \sum_{j \rightsquigarrow k}} (Q_k^T Q_j \eta_{v \, j}^l - \eta_{v \, k}^l) - u_k \times \eta_{v \, k}^l \; , \quad k=1...N \, ,
\end{equation}
again as in \cite{LUCA2}.
\item If $\eta_{\omega \, k}^l \neq 0$, then $\eta_{\omega \, k}^l = Q_k^T \eta_{\omega \, k}^r$ and $\eta_{v \, k}^l = Q_k^T \eta_{v \, k}^r - (Q_k^T r_k) \times (Q_k^T \eta_{\omega \, k}^r)$, and $O_{\C} = \lbrace (\gamma + \beta \times \alpha,\, \alpha) : \alpha, \beta, \gamma \in \mathbb{R}^3 \text{ and } \Vert \gamma \Vert \leq 1 \rbrace$. Designing a consensus algorithm, that achieves agreement on $\xi^r \in O_{\C}$ \emph{and} can be written with left-invariant variables, appears to be difficult. Similarly to the first case, suitable algorithms can be built if the overall dimension of the variables used for the consensus algorithm is enlarged with respect to the dimension of the configuration space. The consensus algorithm proposed in \cite{LUCA2} replaces $\eta_k^l$ by three components $\alpha_k = \eta_{\omega \, k}^l \in \mathbb{R}^3$, $\beta_k \in \mathbb{R}^3$ and $\gamma_k \in \mathbb{R}^3$ associated with the vectors $\alpha$, $\beta$, $\gamma$ used to describe $O_{\C}$ above; then $\eta_k^l = (\eta_{v \, k}^l,\, \eta_{\omega \, k}^l) = (\gamma_k + \beta_k \times \alpha_k,\, \alpha_k)$. The advantage of this embedding $\eta_k^l \rightarrow (\alpha_k,\beta_k,\gamma_k)$ is that left-invariant consensus algorithms can be decoupled for the $\alpha_k$, the $\beta_k$ and the $\gamma_k$. With the notations of the present paper, the corresponding consensus algorithm proposed in \cite{LUCA2} is
\begin{eqnarray*}
\tfrac{d}{dt} \alpha_k & = & {\textstyle \sum_{j \rightsquigarrow k}} (Q_k^T Q_j \alpha_j - \alpha_k) \; - u_k \times \alpha_k\\
\tfrac{d}{dt} \beta_k & = & {\textstyle \sum_{j \rightsquigarrow k}} (Q_k^T Q_j \beta_j - \beta_k + Q_k^T(r_j-r_k)) \; - u_k \times \beta_k - \mathbf{e}_1\\
\tfrac{d}{dt} \gamma_k & = & {\textstyle \sum_{j \rightsquigarrow k}} (Q_k^T Q_j \gamma_j - \gamma_k) \; - u_k \times \gamma_k \phantom{KKKKKKKKKKKK} \; , \quad k=1...N \, .
\end{eqnarray*}
Comparing the left-invariant relative position $g_k^{-1} g_j = (Q_k^T(r_j-r_k),\, Q_k^T Q_j)$ with the terms and factors appearing in this consensus algorithm, one observes that the latter is indeed left-invariant. It can be verified (see \cite{LUCA2}) that this algorithm indeed synchronizes the $\eta_k^r = Ad_{g_k}(\gamma_k + \beta_k \times \alpha_k ,\, \alpha_k)$.\\
\end{itemize}

\noindent \emph{Remark 3:} \LIC in linear motion, i.e. with $\eta_{\omega \, k}^l = 0$ $\forall k$, under steering control requires to align vectors $Q_k \mathbf{e}_1$ for all agents. This is in fact equivalent to \TC on $SO(3)$ with $\eta_k^l = \omega^l = \mathbf{e}_1$ $\forall k$. The present section thus illustrates the method for \TC on $SO(3)$ for uniformly connected $\G$ (instead of fixed undirected $\G$ as in Section 4).\\

\noindent \emph{Remark 4:} \LIC under steering control on $SE(2)$ is treated in \cite{TAC2,SPL2005}. As for $SE(3)$, control algorithms obtained intuitively, with several simplifications due to the lower dimension, can be recovered with the general method of the present paper.

In fact, the group structure and control setting of steering control on $SE(2)$ are such that $\forall g \in SE(2)$ and $\forall$ steering controls $u \in \mathbb{R}$, one has
\begin{equation}\label{zzttyy}
\xi^r = Ad_g \xi^l = Ad_{g}(a+Bu) = \alpha(g,u) + Bu \; \; \text{ with } \; \; \alpha(g,u) \perp Bu \, ;
\end{equation}
On $SE(2)$ explicitly, $a+Bu = (\mathbf{e}_1,u) \in \mathbb{R}^2 \times \mathbb{R}$ and $Ad_g(\mathbf{e}_1,u) = (Q_{\theta} \mathbf{e}_1- u Q_{\pi/2}r,\, u)$, so $\alpha(g,u) = (Q_{\theta} \mathbf{e}_1- u Q_{\pi/2}r,\, 0)$ and $Bu = (0,\, u)$. Then \LIC automatically implies equal $u_k$, thus \RIC, meaning that \emph{underactuated  \LIC is equivalent to \TC} and imposes the same constraints on relative positions $\lambda_{jk}$. This is the case for any group and control setting satisfying (\ref{zzttyy}).

For steering control on $SE(3)$, \LIC is slightly different from \TC because $Ad_g(\mathbf{e}_1,u) = (Q \mathbf{e}_1 + r \times (Q u), Q u)$, so (\ref{zzttyy}) would require $(Q u) \cdot (Q \mathbf{e}_1) = u \cdot \mathbf{e}_1 = 0$ which is not true in general. Therefore, for \LIC under steering control the $\omega^l_k = u_k$ can differ by arbitrary rotations around $\mathbf{e}_1$, while \TC would require equal $\omega^l_k$.



\section{Conclusion}

This paper proposes a geometric framework for coordination on general Lie groups and related methods for the design of controllers driving a swarm of underactuated, simple integrator agents towards coordination. It shows how this general framework provides control laws for coordination of rigid bodies, on groups $SO(3)$, $SE(2)$ and $SE(3)$, and allows to easily handle different settings.

Following the numerous results about coordination on particular Lie groups, various directions are still open to extend the general framework of the present paper. A first case often encountered in practice is to stabilize \emph{specific relative positions} of the agents (``formation control''). In \cite{SPL2005,TAC2} for instance, the steering controlled agents on $SE(2)$ are not only coordinated on a circle, but regular distribution of the agents on the circle is also stabilized; in the present paper, relative positions of the agents are asymptotically fixed but arbitrary. The requirement of synchronization (most prominently, ``attitude synchronization'' on $SO(3)$) also fits in this category. A second important extension would be to consider \emph{more complex dynamics}, like those encountered in mechanical systems.



\bibliographystyle{abbrv}
\bibliography{LGC_Bib.bib}
\end{document}